\documentclass[12pt]{article}
\usepackage{amsmath,amssymb}
\usepackage{color}
\usepackage{theorem}
\usepackage[english]{babel}
\usepackage[latin1]{inputenc}
\usepackage{latexsym}
\usepackage{amscd}
\usepackage{mathrsfs}
\usepackage{indentfirst}
\usepackage[dvips]{graphicx}

\numberwithin{figure}{section}

\newcommand{\R}{\ensuremath{\mathbb{R}}}
\newcommand{\N}{\ensuremath{\mathbb{N}}}
\newcommand{\cs}{contact structure}
\newcommand{\css}{contact structures}

\newcommand{\al}{\ensuremath{\alpha}}
\newcommand{\lt}{\\ \indent}

\newcommand{\ssi}{if and only if}
\newcommand{\SU}{\ensuremath{\mathbb{S}^1}}

\newcommand{\SB}{\ensuremath{\mathscr{C}}}

\newcommand{\carr}{\hfill{$\square$}}
\newcommand{\demo}{\noindent \bf Proof \rm}

\newcommand{\In}{[0,1]}
\newcommand{\cI}{\times [0,1]}
\newcommand{\sbr}{branched surface}
\newcommand{\sbrs}{branched surfaces}
\newcommand{\B}{\ensuremath{\mathscr{B}}}
\newcommand{\Su}{\ensuremath{\mathscr{S}}}

\newcommand{\di}{\ensuremath{\mathscr{D}}}

\newcommand{\Ls}{\ensuremath{\mathscr{L}}}

\newcommand{\ad}{\ensuremath{\mathfrak{a}}}
\newcommand{\NB}{\ensuremath{N(\mathscr{B})}}
\newcommand{\dhNB}{\ensuremath{\partial _h N(\mathscr{B})}}

\newcommand{\dvNB}{\ensuremath{\partial _v N(\mathscr{B})}}
\newcommand{\dvS}{\ensuremath{\partial _v N(\Sigma)}}

\newcommand{\ls}{singular locus}
\newcommand{\vf}{fibred neigbourhood}

\newcommand{\suc}{surface of contact}

\newcommand{\supu}{sink surface}
\newcommand{\suct}{twisted surface of contact}
\newcommand{\dix}{twisted disk of contact}

\newcommand{\sucts}{twisted surfaces of contact}

\newcommand{\dec}{splitting}

\newcommand{\eps}{\ensuremath{\varepsilon}}

\newcommand{\V}{\ensuremath{\mathscr{V}}}
\newcommand{\maxv}{\ensuremath{max(\V)}}
\newcommand{\maxvu}{\ensuremath{max(\V_1)}}
\newcommand{\ac}{\ensuremath{\mathscr{A}}}
\newcommand{\tr}{train track}
\newcommand{\Bs}{\ensuremath{\B_{sup}}}

\newcommand{\Bn}{\ensuremath{\B_{nul}}}
\newcommand{\ga}{\ensuremath{\gamma}}
\newcommand{\gs}{\ensuremath{\gamma_{sup}}}
\newcommand{\dn}{\ensuremath{\delta_{nul}}}
\newcommand{\vad}{\ensuremath{\mathscr{V}(\ad)}}

\newcommand{\twc}{twisted curve}

\newtheorem{df}[figure]{Definition}
\newtheorem{leme}[figure]{Lemma}
\newtheorem{prop}[figure]{Proposition}
\newtheorem{teo}[figure]{Theorem}
\newtheorem{ex}[figure]{Example}
\newtheorem{rmk}[figure]{Remark}
\newtheorem{cor}[figure]{Corollary}

\newcommand{\defi}{\begin{df} \rm}
\newcommand{\rke}{\begin{rmk} \rm}
\newcommand{\expl}{\begin{ex} \rm}

\begin{document}
\renewcommand{\bibname}{References}
\vspace*{1cm}
\Large \bf \begin{center} 
A SUFFICIENT CONDITION FOR A BRANCHED SURFACE TO FULLY CARRY A LAMINATION\vspace{.5cm}\\
\large
\rm
Skander ZANNAD
\end{center}
\normalsize\rm
\vspace{.8cm}
\begin{quote} \bf Abstract - \rm We give a sufficient condition for a \sbr\ in 
a 3 dimensional manifold to fully carry a lamination, giving a piece of answer 
to a classical question of D. Gabai.\\

\noindent \bf Key words : \rm branched surface ; lamination ; twisted disk of
contact ; twisted curve.
\end{quote}

\vspace{.8cm}
\setcounter{section}{-1}
\section{Introduction}\label{ch0}

Branched surfaces are combinatorial objects which prove to be useful, in particular to study laminations. 
They are the main tool to construct essential laminations in the works of D. Gabai, U. Oertel, A. Hatcher, 
T. Li or C. Delman and R. Roberts for instance.\lt 
One of the most striking topological results is theorem \bf \ref{cover} \rm of \bf [GO] \rm :

\begin{teo}\label{cover} \bf ([GO]) \it If a compact orientable 3 dimensional manifold $M$ admits an essential lamination, then its universal cover is homeomorphic to $\R^3$.
\end{teo}

The characterisation of the \sbr s fully carrying an essential lamination is 
now known, after works of 
D. Gabai and U. Oertel, and of T. \hbox{Li :}

\begin{teo}\label{gob} \bf ([GO]) \it A lamination is essential \ssi\ it is fully carried by an essential \sbr.
\end{teo}

\begin{teo} \bf ([Li]) \it  Let $M$ be a closed and orientable manifold. Then every laminar \sbr\ in $M$ fully
carries an essential lamination, and any essential lamination which is not a lamination by planes is 
fully carried by a laminar \sbr.
\end{teo}

D. Gabai and U. Oertel gave a number of necessary and sufficient conditions for a \sbr\ to fully carry an essential lamination, assuming that this \sbr\ already fully carries a lamination. The important contribution of T. Li was to give a sufficient condition for this kind of \sbr\ to fully carry a lamination, highlighting the importance of the following problem of D. Gabai (problems \bf 3.4 \rm of \bf [GO] \rm and \bf 2.1\rm of \bf [Ga]\rm): \it when does a \sbr\ fully carry a lamination?\rm\lt

This question is complex, as shown by L. Mosher's theorem :
\begin{teo}\label{mosh} \bf (L. Mosher) \it The problem of whether or not a general branched surface abstractly carries a lamination is algorithmically unsolvable.
\end{teo}

Let us give brief explanations of the terms ``general \sbr " and ``abstractly carries" : the \sbr s we will use in this text are by definition embedded in a 3-manifold. However, a general \sbr\ could be defined the same way, but without assuming it is embedded or even immersed in a 3-manifold. In \bf[Ch]\rm , J. Christy gives necessary and sufficient conditions for a general \sbr\ to be immersed or embedded in a 3-manifold, and some examples. ``Abstractly carrying" a lamination is the generalization for general \sbr s of ``fully carrying" a lamination. Precise definitions can be found in \bf[MO]\rm . A proof of theorem \bf\ref{mosh} \rm is given in \bf[Ga]\rm.
\lt

The goal of this article is to prove the following result, which is a piece of answer to the question of D. Gabai :

\begin{teo}\label{but}
Let $M$ be an oriented manifold of dimension 3, without boundary. Let \B\ be an orientable \sbr\ of $M$ without twisted curve. Then \B\ fully carries a lamination.\end{teo}

A corollary easily comes from this theorem :

\begin{cor}\label{corobut} Let $M$ be an oriented manifold of dimension 3, without boundary. Let \B\ be an 
orientable \sbr\ of $M$ without twisted curve homotopic to zero in $M$. Then the lift of \B\ in the universal
cover $\tilde{M}$ of $M$
fully carries a lamination\end{cor}

This result is almost optimal in the following sense : the existence of a twisted curve homotopic to zero 
implies the existence of a closed curve homotopic to zero and transverse to \B. But, according to point (4) in lemma \textbf{2.7} of \textbf{[GO]}, if \B\ fully carries an essential lamination, such a closed curve cannot exist. The condition ``there is no twisted curve homotopic to zero" is then sufficient for the lift of \B\ in the universal cover of $M$ to fully carry a lamination, but it is also necessary for this lift to fully carry an essential lamination.\lt

This sufficient condition appeared when investigating on the notion of \cs\ carried by a \sbr. There is some 
hope to use this criterion in the study of \css\ \it via \rm \sbr s.\\

The basic definitions about \sbr s, surfaces of contact and \twc s are given in section \textbf{\ref{ch1}}. 
The principle of the proof of theorem \textbf{\ref{but}} is the same as the one of the construction of a lamination
whose holonomy is strictly negative,
in section \textbf{4} of \bf[OS2]\rm. We will build a resolving sequence
of \dec s, whose inverse limit
induces a null holonomy lamination on the \vf\ of the neighbourhood of
the 1-skeleton of some cell
decomposition into disks and half-planes of \B. Splittings, resolving sequences and inverse limits are introduced in section \textbf{\ref{split}}. Theorem \textbf{\ref{but}} is then proved in section \textbf{\ref{demolem}}.
The last section contains some remarks about the question of D. Gabai.\\

At last, I'd like to thank U. Oertel and J. \'Swiatkowski for having found a mistake in an optimistic version of this text, in which I thought I had answered the question of D. Gabai.

\numberwithin{figure}{subsection}

\section{Branched surfaces}\label{ch1} 

Throughout this article, $M$ is a 3 dimensional oriented manifold without boundary. It will be supposed paracompact (because of remark \textbf{\ref{paracompact}}) and separated.
Its universal cover is denoted $\tilde{M}$.
\subsection{First definitions}

\begin{figure}[!h]
\begin{center}
\input{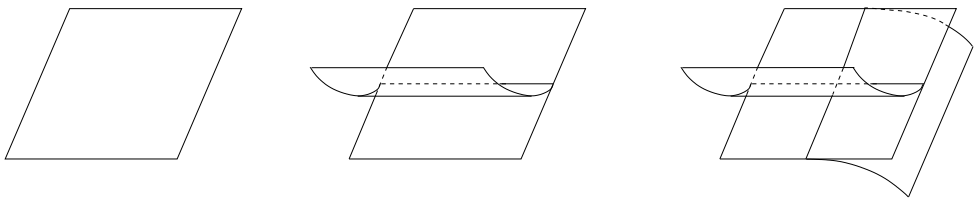}
\caption{Local models of a branched surface}
\label{branch}
\end{center}
\end{figure}

\defi A \textit{\sbr} \B\ in $M$ is a union of smooth surfaces locally
modeled on one of the three models of
figure \textbf{\ref{branch}}. The \textit{\ls}
\Ls\ of \B\ is the set of points, called \textit{branch points}, none of
whose neighbourhoods is a disk.
Its \textit{regular part} is
$\B \backslash \Ls$. The closure of a connected component of the regular
part is called a \textit{sector} of \B.
\end{df}

The \ls\ may have double points, as it is the case in the third model of
figure \textbf{\ref{branch}}. 

\rke\label{paracompact}  According to the local models, the double points are isolated, and since $M$ is paracompact, they are countable.
\end{rmk}

At each regular point of \Ls, we can define a
\textit{branch direction}, as on figure \textbf{\ref{sens}}.

\begin{figure}[!h]
\begin{center}
\input{sens.pstex_t}
\caption{branch direction}
\label{sens}
\end{center}
\end{figure}

\defi\label{voisina} A \textit{\vf} \NB\ of \B\ is an interval ``bundle" over \B,
as seen on figure \textbf{\ref{bord}}. The boundary of
\NB\ can be decomposed into an
\textit{horizontal boundary} \dhNB\ transverse to the fibres and a
\textit{vertical boundary} \dvNB, tangent to the fibres (see figure
\bf\ref{bord}\rm, \textbf{a)}).

\begin{figure}[!h]
\begin{center}
\input{bord.pstex_t}
\caption{Fibred neighbourhood of \B}
\label{bord}
\end{center}
\end{figure}

We define the projection map $\pi\ :\ \NB\to\B$ which sends a fibre of
\NB\ onto its base point.
In particular, $\pi(\dvNB)=\Ls$.
We can also consider \NB, not as an abstract bundle but rather as a part of $M$, and in this case \B\ is not included in \NB.
However, \NB\ contains a \sbr\ $\B_1$ which is isomorphic to \B\ (see
figure \textbf{\ref{bord}}, \textbf{b)}). The \sbr\ $\B_1$ is a \it\dec\ \rm of \B\ (splittings
will be defined in section
\textbf{\ref{split}}).

\end{df}

Let's see how we can put a sign on each double point of \Ls.
Locally, two smooth parts of \Ls\ run through $p$. They are cooriented
by their branch direction, and we
call them $\Ls_1$ and $\Ls_2$.
\begin{figure}[!h]
\begin{center}
\input{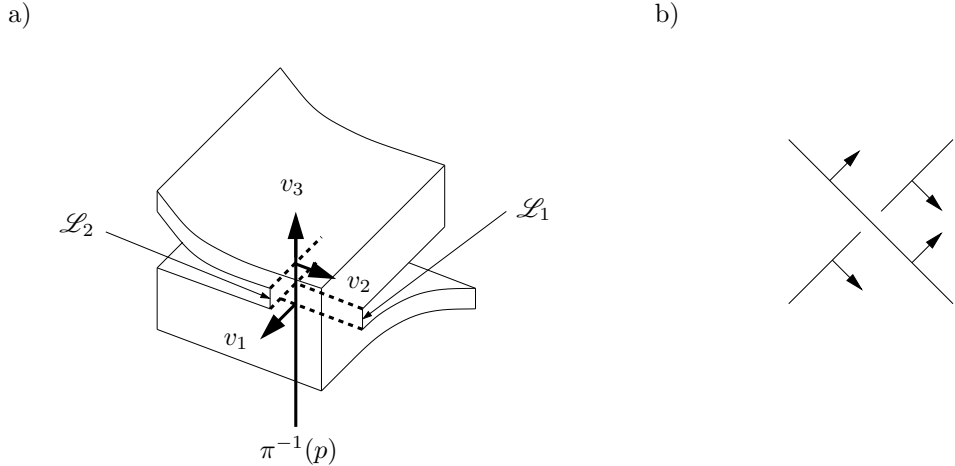}
\caption{a) : $\{v_1,v_2,v_3\}$ ; b) positive double point}
\label{bobo}
\end{center}
\end{figure}
Set an orientation of the fibre of \NB\ passing through $p$.
Hence, it makes sense to say that one of the branching $\Ls_1$ or
$\Ls_2$ is over the other at $p$.
Say for example that $\Ls_1$ is under $\Ls_2$. Let
$v_1$ be a vector of $T_pM$ defining the branch direction
of $\Ls_1$ at $p$, and $v_2$ be a vector of $T_pM$ defining the
branch direction
of $\Ls_2$ at $p$. At last, let $v_3$ be a vector giving the chosen
orientation of the fibre of \NB\ passing
through $p$, as seen on figure \bf\ref{bobo}\rm, \textbf{a)}. We then
call $p$ a
\it positive double point
\rm (resp. \it negative double point\rm) if the base $\{v_1,v_2,v_3\}$
of $T_pM$ is
direct (resp. indirect) in respect with the orientation of $M$.
With this convention, the positive double points will be
drawn in the plane as on the diagram
\textbf{b)} of figure \bf\ref{bobo}\rm.

\rke The sign of a double point depends on the orientation of $M$: if
this one is reversed, the signs of the
double points are reversed as well.
Though, this sign is independent of the chosen orientation of the fibre
passing through the double point in
the preceding definition.
\end{rmk}

\defi A codimension 1 \it lamination \rm in a dimension 3 (resp. 2)
manifold $M$ is the decomposition of a
closed subset $\lambda$ of $M$ into injectively immersed surfaces (resp. curves) called \it
leaves\rm,
such that $\lambda$ is covered
by charts of the form $]0,1[^2\times I$ (resp. $]0,1[\times I$) in which the
leaves have the form
$]0,1[^2\times \{point\}$ (resp. $(]0,1[\times \{point\}$).
\end{df}

\defi A \sbr\ \B\ \it carries a lamination \rm $\lambda$ of codimension
1 if $\lambda$ is contained in a
\vf\ of \B\ and if its leaves are transverse to the fibres. We say that
$\lambda$ is \it fully carried \rm
if moreover it meets all the fibres.
\end{df}

\subsection{Surfaces of contact}\label{sucs}

Let \B\ be a \sbr.

\defi A \it\suc\ \rm\ is the immersion of a surface \Su\ in \B, whose boundary 
is sent onto smooth circles of the \ls\ of \B, such that the branch directions along these boundary components
point into \Su.\lt
If we consider a lift of \Su\ into \NB, we see that the existence of such a surface is equivalent to the
existence of
an immersion \hbox{$f:\Su\to\NB$} satisfying :
\begin{itemize}
\item[(i)] $f(Int(\Su))\subset Int(\NB)$ and is transverse to the fibres ;
\item[(ii)] $f(\partial{\Su})\subset Int(\dvNB)$ and is
transverse to the fibres.\vspace{.2cm}\lt
\end{itemize}
Hence, the expression \it \suc\ \rm will be used for both definitions.

\end{df}

\begin{figure}[!h]\begin{center}
\input{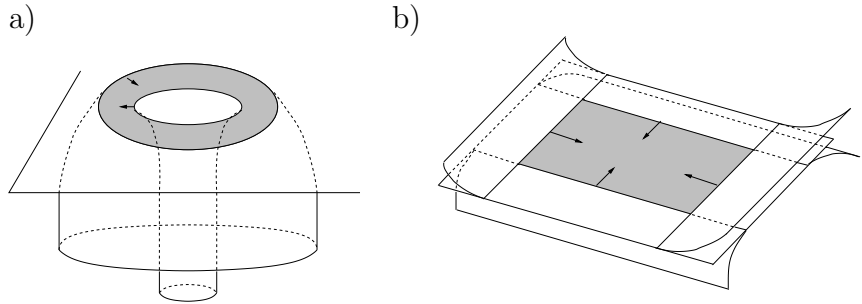}\caption{Annulus of contact and sink
disk}\label{isc}\end{center}\end{figure}

An example is given in figure \bf\ref{isc}\rm, \textbf{a)}.

\rke In general, a surface of contact is not a sector, but a union of
sectors : the \ls\ of the \sbr\
may meet the interior of the surface of contact.
The same is true for the sink surfaces and the \sucts\ defined further.
\end{rmk}

\defi A \it\supu\rm\ is the immersion of a surface \Su\ in \B, whose boundary 
is sent onto piecewise smooth circles of the \ls\ of \B, at least one of whose is not smooth, such that the branch directions along these boundary components
point into \Su. A double point in the boundary of \Su\ which is the intersection of two smooth components of
the boundary of \Su\ is called a \it corner \rm of \Su. A \supu\ has thus at least one corner.\lt
Equivalently, if we consider a non smooth lift of \Su\ into \NB, we can say that a \supu\ is
an immersion \hbox{$f:\Su\to\NB$} satisfying :

\begin{itemize}
\item[(i)] $f( Int(\Su))\subset Int(\NB)$ and is transverse to the fibres ;
\item[(ii)] $f(\partial{\Su})$ is included in $Int(\dvNB)$
except in a finite and non empty number of closed intervals $C_1\ldots C_k$.
Outside these $C_i$,
$f(\partial{\Su})$ is transverse to the fibres of
$\dvNB$.
Each $C_i$ is included in a fibre of \NB\ corresponding to a double point of
\Ls, and must intersect $Int(\NB)$. Thus,
$\pi(f(\partial{\Su}))$ is not smooth.
The $C_i$ s are called the
\it corners \rm of \Su.
\vspace{.2cm}
\end{itemize}
\end{df}

An example is given in figure \bf\ref{isc}\rm, \textbf{b)}.

 \defi A \textit{\suct} is a \supu\ whose every corners, which are double
points, have the same sign, and which satisfies : for some Riemannian metric for which $\Ls$ in the neighbourhood of a double point cuts \B\ into four
sectors of angle $\pi/2$, the corners of a \suct\ are all of angle $\pi/2$. The case of a corner of angle $3\pi/2$ is forbidden.\lt
A \suct\ is \it positive \rm (resp. \it negative\rm ) if all its corners are
positive (resp. negative).
\end{df}

An example is given in  figure \bf\ref{isc2}\rm.

\vspace{.3cm}
\begin{figure}[!h]\begin{center}\input{isc2.pstex_t}\caption{Negative \dix }
\label{isc2}\end{center}\end{figure}

\rke\label{nonplong} It is well-known that the existence of a \dix\ which is a sector is an obstruction to the existence of a lamination fully carried. This fact is stated in proposition \textbf{\ref{nonport}}, and then proved, by using \textit{\tr s}.\end{rmk}

\subsection{Twisted curves}

\defi A \sbr\ \B\ is said \textit{orientable} if there exists a global orientation of the fibres of a \vf\ \NB\ of \B.
\end{df}

\rke An orientable \sbr\ cannot have any \textit{monogon}, i.e. a disk $D\subset M\backslash\ Int(\NB)$ with $\partial D=D\cap\NB=\beta\cup\delta$, where $\beta\subset\dvNB$ lies in a fibre of
\dvNB\ and 
$\delta\subset\dhNB$ (see figure \bf\ref{monog}\rm ).

\begin{figure}[!h]\begin{center}\input{monog.pstex_t}\caption{monogon}\label{monog}
\end{center}\end{figure}
\end{rmk}

Throughout this section, \B\ will now be an orientable \sbr, and \NB\ a \vf\ of \B, with a fixed orientation of the fibres.

\defi A \textit{positive} (resp. \textit{negative})  \textit{twisted curve} $\gamma$ is an oriented closed curve, immersed in \B\ and which satisfies :

\begin{itemize}
\item [(i)] $\gamma$ is included in \Ls , and is then a finite union of smooth segments of \Ls. When, along \ga , we pass from a smooth segment to another one through a double point of \Ls, this double point is called a \it corner \rm of \ga\ ;
\item[(ii)] $\gamma$ has at least one corner ;
\item[(iii)] at a corner, \ga\ passes from a smooth segment $I_1$ of \Ls\ to a smooth segment $I_2$ of \Ls. Since \B\ is oriented, one of these segments is over the other one. If $I_2$ is over $I_1$, we then say that the corner is \it ascending\rm , else it is \it descending\rm . The sign of a corner as a double point does not determine whether it is ascending or descending. We then demand all the corners of \ga\ to be ascending (resp. descending).
\end{itemize}
The existence of a twisted curve \ga\ is equivalent to the existence in \NB\ of a curve, still denoted \ga , which can be decomposed into a union $\gamma=\gamma_l\cup\gamma_c$, where $\gamma_l$ and $\gamma_c$ verify :
\begin{itemize}
\item[(i)] $\gamma_l$ is the \textit{smooth part} of \ga\ : it is a finite union of segments included in \dvNB\ and transverse to the fibres of \NB\ ; 
\item[(ii)] $\gamma_c$ is a finite and non empty union of segments denoted $C_i,\ i=1\ldots n$, where each $C_i$ is included in a fibre of \NB, in such a way that the orientation of $C_i$, induced by the one of \ga , coincides with (resp. is opposite to) the orientation of this fibre. The $C_i$ 's are the corners of \ga. They are said ascending (resp. descending) if \ga\ is positive (resp. negative).
\end{itemize}
\end{df}

\rke If we revert the orientation of a positive \twc, we get a negative \twc. The converse is also true.
\end{rmk}

\rke The boundary of a \suct\ is a \twc.\end{rmk}

\defi A \it simple \rm \twc\ is a \twc\ whose inside of the smooth part is embedded in \B. Otherwise said, 
only the corners are of multiplicity 2 or more.\end{df}

\begin{leme}\label{simplet} Let \ga\ be a \twc . Then there exists a simple \twc\ $\delta$ included in \ga.
\end{leme}

\demo We also denote $\gamma\ : \SU\to\B$ the immersion of the \twc\ \ga , \SU\ being oriented. We can always suppose that \ga\ is positive. If \ga\ is simple, we obviously have $\delta=\ga$. 
Else, there exist two segments of \SU\ denoted $J=[a,b]$ and $K=[c,d]$, whose interiors are disjoint, and whose orientation is the one induced by the orientation of \SU, and such that the images $\ga(J)$ and $\ga(K)$ coincide.
If $\ga(J)=\ga(\SU)$, we then set $\delta_1=\ga(J)$, which is a \twc.\lt
Else, the segments $J$ and $K$ are chosen to be maximal, in the sense that for every sufficiently small neighbourhood $\mathcal{V}(J)$ of $J$ in \SU\ and every sufficiently small neighbourhood $\mathcal{V}(K)$ of $K$ in \SU, we have
$\ga(\mathcal{V}(J))\not\subset\ga(\mathcal{V}(K))$ and $\ga(\mathcal{V}(K))\not\subset\ga(\mathcal{V}(J))$. 
This means that $\ga(\mathcal{V}(J))$ and $\ga(\mathcal{V}(K))$ split at $\ga(a)$ and at $\ga(b)$. Thus, the point $A=\ga(a)$ of \Ls\ is a corner of $\ga(\mathcal{V}(J))$ or of $\ga(\mathcal{V}(K))$, and it is the same for 
$B=\ga(b)$.\lt

We then distinguish two cases : \lt

\begin{itemize}
\item[(i)] The orientations of $\ga(J)$ and $\ga(K)$ coincide, that is $\ga(c)=\ga(a)=A$ and $\ga(d)=\ga(b)=B$ :

\end{itemize}

The point $B$ is a corner of $\ga(\mathcal{V}(J))$ or of $\ga(\mathcal{V}(K))$. For example, let us suppose it is a corner of
$\ga(\mathcal{V}(J))$. We then set : $\delta_1=\ga([a,c])$. Since $\ga(c)=\ga(a)$, this curve is closed. Since $c\notin Int(J)$, we have $J=[a,b]\subset [a,c]$, and hence $B$ is a corner of $\delta_1$, because it is a corner of $\ga(\mathcal{V}(J))$. At last, all the corners of $\delta_1$ are ascending since $[a,c]$ is oriented after \SU, and all the corners of \ga\ are ascending. Therefore, $\delta_1$ is a \twc.\lt

\begin{itemize}
\item[(ii)] The orientations of $\ga(J)$ and $\ga(K)$ are opposite, that is $\ga(d)=\ga(a)=A$ and $\ga(c)=\ga(b)=B$ :
\end{itemize}

We then set : $\delta_1=\ga([c,b])$. Since $\ga(c)=\ga(b)$, this curve is closed. Since $B$ is a corner of $\ga(\mathcal{V}(J))$ or of $\ga(\mathcal{V}(K))$, $B$ is a corner of $\delta_1$. At last, all the corners of $\delta_1$ are ascending since $[c,b]$ is oriented after \SU, and all the corners of \ga\ are ascending. Therefore, $\delta_1$ is a \twc.\lt

In every case, \ga\ contains a positive \twc\ $\delta_1$ which has strictly less corners (counted with multiplicity) than \ga. If $\delta_1$ is not simple, we iterate the previous construction to $\delta_1$, and we get a positive \twc\  $\delta_2$, having strictly less corners than $\delta_1$. Since \ga\ is closed, it is compact and has a finite number of corners, even with multiplicity. Therefore, in a finite number of steps, we get a positive simple \twc\ included in \ga.\vspace{.5cm}
\carr

The following corollary follows easily :
\begin{cor}\label{simple} A \sbr\ without simple \twc\ does not have any \twc\ at all.
\end{cor}

\section{Splittings}\label{split}
\subsection{Definitions}

\defi\label{deco} Let \B\ and \B$'$ be two \sbrs\ in $M$. We say that
\B' is a \textit{\dec} of
\B\ if there exists a \vf\ \NB\ of \B\ and
an $I$-bundle $J$ in \NB, over a subsurface of \B, such that :
\begin{itemize}
\item[(i)] $\NB=N(\B')\cup J$ ;
\item[(ii)] $J\cap N(\B')\subset\partial J$ ;
\item[(iii)] $\partial _h J\subset\partial _h N(\B')$ ;
\item[(iv)] $\partial _v J\cap N(\B')$ is included in $\partial _v
N(\B')$, has finitely many components,
and their fibres are fibres of $\partial _v N(\B')$.
\end{itemize}

\end{df}

\begin{rmk}\rm When \B$'$ is a \dec\ of \B, the following notation will be
used : $\B'\overset{p}{\to}\B$.
Actually,
\B$'$ is included in a \vf\ \NB\ of \B, endowed with a projection $\pi$
on \B, and the restriction $p$ of
$\pi$
to \B$'$ is the projection we wanted.
\end{rmk}

\begin{figure}[!h]\begin{center}\input{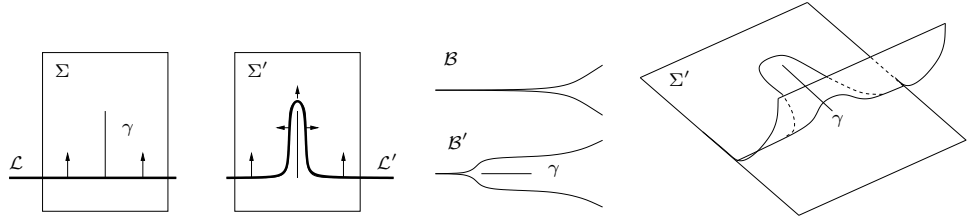}\caption{Splitting
along
$\gamma$}\label{degamma}\end{center}\end{figure}

\defi\label{gam} Let \B\ be a \sbr. Let $\Sigma$ be a sector of \B\
whose boundary contains a smooth part \al\
of
\Ls\ and whose branching direction points into $\Sigma$.
Let $\gamma : I\to\Sigma$ be an embedded arc in
$\Sigma$ such that $\gamma (0)\in\al$ and
\hbox{$\gamma(t)\in Int(\Sigma)$} for $t\neq 0$.
A \textit{\dec\ along $\gamma$} is a \sbr\ \B$'$ defined as in definition
\bf\ref{deco}\rm, where $J$ is an $I$-bundle over a tubular neighbourhood
of $\gamma$ in $\Sigma$ (see
figure \bf\ref{degamma}\rm).

\end{df}

\begin{figure}[!h]\begin{center}\input{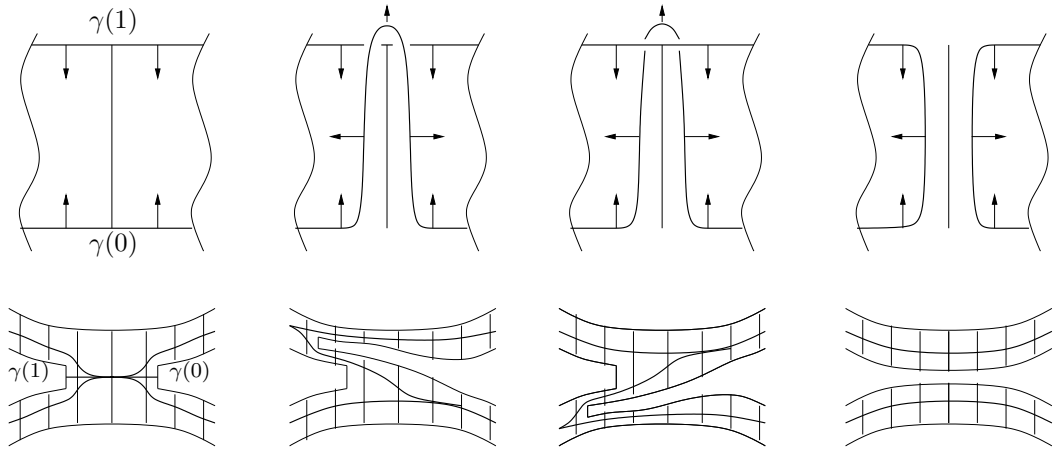}\caption{Over, under
and neutral splittings}\label{3dec}
\end{center}\end{figure}

\defi\label{gam2} We keep the notations of definition \bf\ref{gam}\rm.
Suppose now that
$\gamma(1)$ is in $\Ls$, in a point where the branching direction points
into $\Sigma$ as well.

We then say that $\gamma$ is in \it face-to-face \rm position. If an
orientation of the fibres of \NB\ along
$\gamma$ is chosen, there are three possible splittings along $\gamma$ :
the \it over splitting\rm,
the \it under splitting\rm\ and the \it neutral splitting\rm, drawn in figure
\bf\ref{3dec}\rm.
\end{df}

\rke If $\Sigma$ is a non compact sector, a splitting can be performed
along a non compact arc
$\gamma\ :\  [0,1[\to \Sigma $, verifying the same conditions
as in definition \bf\ref{gam}\rm.
This splitting can be seen as a neutral splitting ``at infinity".
\end{rmk}

\begin{figure}[!h]\begin{center}\input{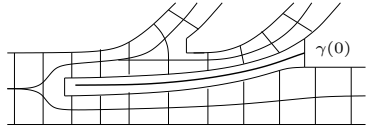}\caption{Backward
splitting}\label{retro}
\end{center}\end{figure}

\rke It is possible to perform a splitting along an arc $\gamma$ which
comes from a sector to another one
through the \ls\ in
the branch direction. In this case, there is only one possible
splitting, called a
\it backward splitting\rm\
(see figure \bf\ref{retro}\rm). 
\end{rmk}

\subsection{Inverse limit of a sequence of splittings}\label{liminv}

\defi
Let \B\ be a \sbr. A \it sequence of splittings \rm is a sequence
\hbox{$\ldots\ \B_{k+1}\overset{p_{k+1}}{\to}\B_k
\overset{p_{k}}{\to}\ldots\overset{p_{2}}{\to}\B_1\overset{p_{1}}{\to}\B=\B_0$}
of \sbr s
$(\B_i)_{i\in\N}$ such that :
\begin{itemize}
\item[(i)] for all $i$, $\B_{i+1}$ is a splitting of $\B_i$ ;
\item[(ii)] for all $i$, $\B_i$ is endowed with a \vf\ $N(\B_i)$, and
those \vf s are such that
$N(\B_{i+1})$ is contained in $N(\B_i)$.
\end{itemize}
Thus, the fibres of each $N(\B_i)$ are tangent to the fibres of
\NB.
\end{df}

For such a sequence, we denote, for all $k\geq 1$ :
$$P_k=p_1\circ p_2\circ\ldots\circ p_k=\pi|_{B_k}:\B_k\to\B$$
the projection from \B$_k$ onto \B.\lt
We will also denote by $\pi_n\ : N(\B_n)\to\B_n$ the projection along the fibres from $N(\B_n)$ to 
$\B_n$.\lt

The following definition is inspired by \bf[MO] \rm :

\defi\label{resolv} A sequence of \dec s
\hbox{$\ldots\ \B_{k+1}\overset{p_{k+1}}{\to}\ldots\overset{p_{1}}{\to}\B=\B_0$}
is said
\it resolving \rm if it satisfies :
\begin{itemize}
\item[(i)] there exist points of \B\ denoted $(x_i)_{i\in\N}$, a real
number $\rho>0$
and disks embedded in
\B\ denoted
$(d_i)_{i\in\N}$, centred at $x_i$ and of radius $\rho$ for some metric
on \B, such that the $d_i$ 's
cover \B\ ;
\item[(ii)] for all $i\in\N$, there exists a subsequence
$(\B_{\varphi_i(n)})_{n\in\N}$ such that the
branch loci of the \sbr s of this subsequence do not intersect
$\pi^{-1}(d_i)$. That is, for all $k$,
$\B_{\varphi_i(k)}$ does not have any branching over $d_i$ :
the branch points over $d_i$ have been
\textit{resolved}, and $P^{-1}_{\varphi_i(k)}(d_i)$ is thus a union of
disjoint disks.
\end{itemize}
When such a sequence exists, we say that \B\ \textit{admits a resolving
sequence of \dec s}.
\end{df}

\rke In particular, a \sbr\ admitting a resolving sequence of splittings
is \textit{fully splittable} in the
sense of
\bf[GO]\rm.
\end{rmk}

\begin{leme}\label{gomo}\rm\bf ([GO],[MO]) \rm \it
Let \B\ be a \sbr\ admitting a resolving sequence of \dec s. Then \B\
fully carries a lamination.
\end{leme}

\demo : It can be found in p. 84-85 of \bf[MO]\rm. \\

 Let
\hbox{$\ldots\ \B_{k+1}\overset{p_{k+1}}{\to}\ldots\overset{p_{1}}{\to}\B=\B_0$}
be a resolving sequence.\\
Let us define \hbox{$\lambda=\bigcap_{n\in\N}N(\B_n)$.} As an intersection of
closed subsets, $\lambda$
is closed. We will now find an
adapted atlas, whose charts will be the
$\pi^{-1}(d_i)$ 's, where the $d_i$ 's
are the disks from point \textbf{(i)} of definition \bf\ref{resolv}\rm.\\
Let $i\in\N$, and $y\in d_i$. Then $\lambda\cap\pi^{-1}(y)$ is some
closed subset $T$ in
$[0,1]$.
The sequence of splittings being resolving, let us consider the subsequence
$(\B_{\varphi_i(n)})_{n\in\N}$ from point \textbf{(ii)} of definition
\bf\ref{resolv}\rm. Since the
$N(\B_n)$ form a decreasing sequence of closed subsets, we get :
$\lambda=\bigcap_{n\in\N}N(\B_{\varphi_i(n)})$.
But, for all $y$ in $d_i$ and for all integer $n$,
$P^{-1}_{\varphi_i(n)}(x)=P^{-1}_{\varphi_i(n)}(y)$, according to
point \textbf{(ii)} of definition
\bf\ref{resolv}\rm. Hence, for all $i$,
$\lambda\cap\pi^{-1}(d_i)$ is topologically a product $d_i\times T$. If the transversal $T$ contains an
interval $I_T$ whose interior is non empty, we remove $Int(I_T)$ from $T$. We then reduce $T$ to a
transversal $T'=T\backslash Int(T)=\partial T$, whose interior is empty, and which is totally discontinuous.
Hence $\pi^{-1}(d_i)$ is a
laminated chart,
the leaves being the
$\{t\}\times d_i$ s, for $t\in T'$. The set $\lambda '=\cup_{i\in\N}(d_i\times T')$ is a lamination.
Moreover, $\lambda'$ meets all the fibres of \NB\ transversally.
\carr

\defi\label{liminv} Let
\hbox{$\ldots\ \B_{k+1}\overset{p_{k+1}}{\to}\ldots\overset{p_{1}}{\to}\B=\B_0$}
be a resolving sequence of
\dec s. The fully carried lamination $\lambda '=\cup_{i\in\N}(d_i\times T')$ defined in the previous proof
is called the
\it inverse limit \rm of this sequence of \dec s.
\end{df}

\section{Proof of theorem \textbf{\ref{but}} }\label{demolem}

Let \B\ be a \sbr\ satisfying the hypotheses of theorem \textbf{\ref{but}}.

\subsection{Cell decomposition of \B }\label{raff}

The \ls\ \Ls\ cuts \B\ into sectors.
This is a first cell decomposition $X$ of \B. The 2-cells are the
sectors, and are not disks or half-planes
in general. The edges are the smooth parts of \Ls\ having no double
points in their interior and such that :
if an edge is compact, both ends are double points (they may be the same
double point) ; if an edge is
diffeomorphic to $[0,1[$, then its end is a double point ; if an edge is
diffeomorphic to \R, it does not meet
any double point. \lt

But this first decomposition is not fine enough. For
reasons which should become clear
after the statement of lemma \textbf{\ref{voi1}}, this decomposition
must be refined into a decomposition
whose compact cells are disks. This is made by adding as many vertices and 
edges (compact or not) as
necessary. We also add vertices and edges so that the
non compact cells are half-planes, and vertices so that no edge is a loop (i.e. its two 
ends coincide). We denote $Y$ the obtained decomposition.

\rke The ``boundary of a 2-cell" is not the topological boundary, but the
combinatorial one. An edge can be
found twice, with different orientations, in the boundary of the same
2-cell.
\end{rmk}

\subsection{First \dec}\label{prem}

The first step is to perform a first splitting of \B, denoted $\B_1$, which is fully carried by \NB, as in
definition \textbf{\ref{voisina}}. Let us describe it more precisely.
\lt

Let $\eps$ be a non negative real number, such that, for some metric on
\B, the edges of $Y$ are all strictly
longer than $5\eps$ (we shall see why later). Let us look at the
intersection of \B\ with an \eps -tubular
neighbourhood
of \Ls\ in $M$. We chose \eps\ small enough for this tubular neigborhood
to be regular.
This intersection is the union of \Ls\ and of two other parts, which
meet together at the double points :
one part lies behind \Ls, for the coorientation of \Ls\ given by the
branch directions,
and the other part, denoted
$T_\Ls$, lies in front of \Ls.
The boundary of
$T_\Ls$ is included in the union of \Ls\ with a parallel copy of \Ls,
called $\Ls_1$. It is just ``included
in" and not ``equal to" this union, because of what happens at the double
points.
The first \dec\ is a \dec\ over $T_\Ls$, which means that we remove from
\NB\ an $I$-bundle over
$T_\Ls$. The \sbr\ $\B_1$ we get is isomorphic to \B, and its \ls\ is
$\Ls_1$
(see figure
\bf\ref{prems}\rm).

\begin{figure}[!h]\begin{center}\input{prems.pstex_t}\caption{First
splitting}\label{prems}
\end{center}\end{figure}

The trace of $\pi^{-1}(\Ls)$ on \B$_1$ is made of two copies of \Ls,
drawn on two different sectors (at
least locally), as seen on figure \bf\ref{prems}\rm. The trace of
 $\pi^{-1}(Y)$ on \B$_1$,
denoted $Y_1$, is then more complicated than a cell decomposition into
disks and half-planes, since some of
the cells and some of the edges are branched. But all these branchings lie in a closed \eps
-neighbourhood of $\Ls_1$, and \B$_1$ minus a closed
\eps -neighbourhood of the 1-skeleton of $Y_1$ is the same union of disks
and half-planes as \B\ minus the 1-skeleton of $Y$.

\subsection{Train tracks}\label{vf}

Each 2-cell $\Sigma$ of $Y$
inherits from \NB\ an interval bundle
$N(\Sigma)$ built in the following way : we denote by $N(Int(\Sigma))$ the set of all the fibres of \NB\
whose base point lies in $Int(\Sigma)$ and we set 
$N(\Sigma)=\overline{N(Int(\Sigma))}$.

 The boundary of this $N(\Sigma)$
can be decomposed into an horizontal boundary
(included in \dhNB) and a vertical boundary
(not included in \dvNB). 
Since all the compact 2-cells of $Y$ are disks
and are orientable, the vertical
boundary of
$N(\Sigma)$, denoted \dvS, is in fact of the form
$\SU \times I$. For the non compact 2-cells, the vertical boundary is of
the form
$\R\cI$. 
For each 2-cell $\Sigma$, let us look at the trace of $\B_1$ on
\dvS, which is also the boundary of $\Sigma_1=P_1^{-1}(\Sigma)$. It is a
\it train track\rm , i.e. a branched curve fully carried by \dvS.
This \tr\ does not have a boundary and avoids the trace of \dvNB\ on
\dvS. It is compact \ssi\ $\Sigma$ is
compact. Figure \bf\ref{train}\rm\ shows two examples of compact
\tr s.\lt

\begin{figure}[!h]\begin{center}\input{train.pstex_t}\caption{Train
tracks}\label{train}\end{center}
\end{figure}

An orientation of a 2-cell gives an orientation of its boundary. The
corresponding \tr\ is then oriented as
well. For each 2-cell $\Sigma$, we set an orientation of the fibres of
\dvS. We introduce the following
definitions :

\defi A branching of a \tr\ is said \it direct \rm when a track followed
in the direct way divides itself
into two tracks at this branching, and it is said \it backward \rm when
two tracks followed in the direct way
meet at this branching.
\end{df}

We can go a bit further in the classification of the branchings of a \tr\ :

\defi Let \V\ be an oriented compact \tr\ without boundary, fully
carried by a trivial bundle $\SU\cI$. We set an orientation of the fibres. Let \SB\ be a
smooth closed curve of \V. It cuts
$\SU\cI$ into two parts : $(\SU\cI)^+$, containing the points which lie
over \SB\ for the orientation of the
fibres, and $(\SU\cI)^-$ containing the points which lie under. A
branching along a smooth closed curve of
\SB\ is called an \it over branching \rm (resp. \it under branching\rm)
if the branch which leaves or meets
\SB\ there lies in $(\SU\cI)^+$ (resp. $(\SU\cI)^-$).
\end{df}

We can thus state the following lemma :

\begin{leme}\label{voi1}  Let $\Sigma$ be a compact 2-cell of 
$Y$ ($\Sigma$ is a disk), and $\Sigma_1$ be its
trace on $\B_1$. Let \V\ be the \tr\ associated to the boundary of
$\Sigma_1$. It is an oriented compact
\tr\  without boundary fully carried by a bundle $\SU\cI$. We set an orientation of the
fibres.
The three following assertions are equivalent :
\begin{itemize}

\item[(i)] when we follow a smooth closed curve of \V, either no
under branching is met or at least one
direct under branching and one backward under branching are met ;

\item[(ii)] \V\ can be split into a union of smooth circles ;
\item[(iii)] $\Sigma$ is not a \dix.
\end{itemize}
\end{leme}

\demo

\begin{itemize}

\item[$\bigstar\ (i)\Rightarrow (ii)$ :] 
If this is true for each connected
component of \V, then it is true for
\V. So we assume \V\ connected and different from a smooth curve. For $\theta\in\SU$ we define
\hbox{$max(\theta)=$}\hbox{$max\{t\in\In\ |\ (\theta,t)\in\V\}$}, which is in
\In, and then we define
$max(\V)=\{(\theta,max(\theta)),\theta\in\SU\}$. This \maxv\ is a smooth
circle of \V, along which we meet at
least one direct under branching and one backward under branching, and no over
branching. In particular, there exists
an oriented arc \ac\ of \maxv, going (for the orientation of \V), from a
direct branching to a backward
branching, with no branching between the two previous ones. 
Then $\V\backslash\ac$ is an oriented compact \tr\ without boundary
denoted \V$_1$, fully carried by
$\SU\cI$, and $\V=\V_1\cup\maxv$. Each smooth closed curve of \V$_1$ is
a smooth closed curve of \V\ as well,
and its under branchings remain unchanged by the previous splitting.
Hence, \V$_1$ satisfies point
\textbf{(i)} of the lemma. If \V$_1$ is not a circle, we perform the
same operation again using \maxvu, and
after a finite number of steps, we have decomposed \V\ into a union of
smooth circles. An example is shown
in figure \bf\ref{baba}\rm.

\begin{figure}[!h]\begin{center}\input{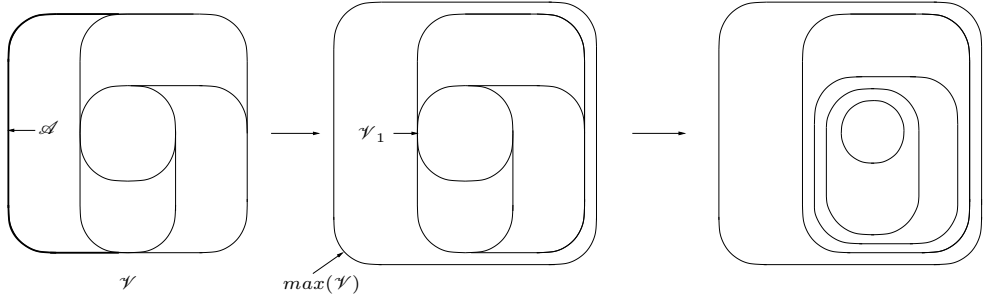}
\caption{Splitting of a \tr\ into a union of smooth
circles}\label{baba}\end{center}\end{figure}

\item[$\bigstar\ \neg (i)\Rightarrow \neg (ii)$ :] 
Let \SB\ be a smooth
closed curve of \V\ having, for
example, only direct under branchings. If we follow \SB\ in the direct
way, and if we take a direct under
branching, then, whatever the smooth path we follow on \V, we will never
be able to go on \SB\ again, for it
would imply the existence of a backward under branching along \SB. Thus,
no branch leaving \SB\ by a direct
under branching is included in smooth closed circle of \V, and \V\ is not a union of smooth circles.

\end{itemize}

\begin{itemize}
\item[$\bigstar\ \neg (iii)\Rightarrow \neg (i)$ :] 
The 
trace of $\B_1$ when $\Sigma$ is a \dix\ is always as on figure 
\bf\ref{bibi}\rm, i.e. it is the union of two smooth circles
and of segments joigning them at branch points.

\begin{figure}[!h]\begin{center}\input{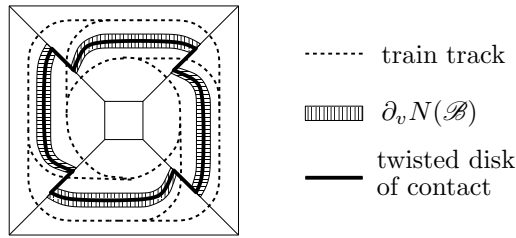}
\caption{$\neg (i)\Leftrightarrow \neg
(ii)$}\label{bibi}\end{center}\end{figure}

The top circle has only under branchings, and it has at least
one under branching because a \dix\ has at least one corner. Moreover, these branchings are of one type because the
corners of a \dix\ all have the same sign.

\item[$\bigstar\ \neg (i)\Rightarrow \neg (iii)$ :] 

Suppose that there
exists a closed smooth curve $C$ of \V\
whose under branchings are all of a single type, for example direct, and which has at least one
under branching. Look at the trace of \dvNB\ on
$\partial_v (\Sigma)$. Each of its connected components has a vertical boundary with two connected components and a
horizontal boundary, also with two connected components. 
Each component of the vertical boundary is included in
a fibre of $\partial_v N(\Sigma)$ whose base point is a branch point 
of \B.\lt
Stand at a point $p$ on $C$, and follow $C$ in the direct way. 
When we meet the first branch point $p_1$, $C$
divides into two branches : the top branch passes over a component $b_1$ of \dvNB, and the bottom branch
passes under $b_1$. We go on until we meet the fibre where $b_1$ ends, whose base point is some branch point
$p_2$. If $p_2$ is not a double point, then the branch of \V\ which is over $b_1$ joins the branch which is
under $b_1$. But these two branches are the two previous branches, and that would imply that there is a
backward branching on $C$. Hence $p_2$ is a double point. At $p_2$, there are thus two branchings, one is 
direct and the other is backwards. One of them is on $C$, so this is the direct one. Again, $C$ divides into two branches which surrounds another component $b_2$ of the trace of \dvNB.
Since this branching is direct, $b_2$ lies over $b_1$ at $p_2$.
We carry on
following $C$ until we return at $p_1$. We have then met $k$ components $b_1\ldots b_k$ of \dvNB\ and $k$ double
points $p_1\ldots p_k$. Each $b_i$ goes from $p_i$ to $p_{i+1}$ for $i=1\ldots k$ modulo $k$. At $p_i$, 
$b_{i-1}$ lies under $b_i$, for all $i$. As a result, 
all the double points have the same sign, and $\Sigma$
is a \dix\ with $k$ corners.
\end{itemize}
\carr

\rke In the proof of point \textbf{(i) $\Rightarrow$ (ii)}, we could also define $min(\V)$ in the same way 
as $max(\V)$, and
show that points \textbf{(ii)} and \textbf{(iii)} are
equivalent to a point \textbf{(i')} : when we follow a smooth closed
curve of \V, either we meet no over
branching, or we meet at least one direct over branching and at least
one backward over branching.
Points \textbf{(i)} and \textbf{(i')} are thus equivalent.
\end{rmk}

With the same ideas, we can also prove the following well known proposition, whose result has already been 
mentioned in remark \textbf{\ref{nonplong}} :

\begin{prop}\label{nonport} Let \B\ be a \sbr\ having a disk sector \di\ which is a \dix\ as well. 
Then \B\ cannot fully carry a lamination.
\end{prop}

\begin{figure}[!h]\begin{center}\input{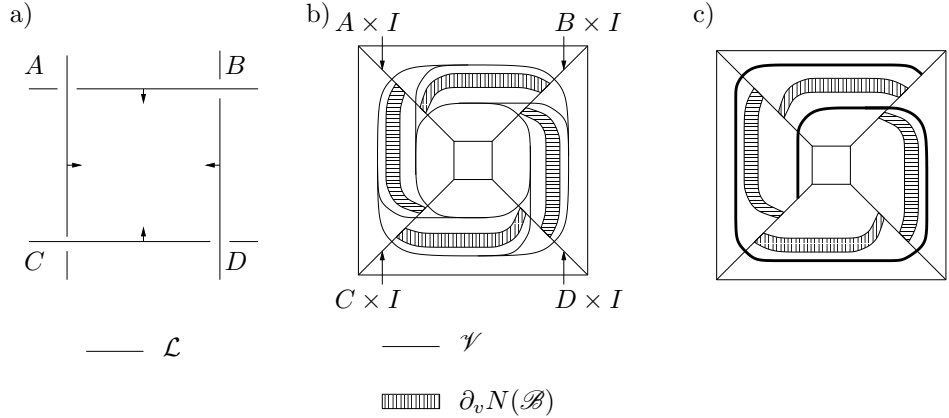}\caption{Twisted
disk of contact}
\label{tordu}\end{center}\end{figure}

\demo
Suppose that \B\ fully carries a lamination $\lambda$. We consider $N(\di)$, the \vf\ over $\di$ in \NB. 
Let \V\ be the \tr\ which is the trace of $\B_1$ in $\partial_v(N(\di))$. The intersection of the 
leaves of $\lambda$ passing over \di\ with $\partial_vN(\di)$ is a union of disks. Their boundaries are circles which form a 
1-dimensional lamination fully carried by some \vf\ of \V. However, as seen in the example on figure 
\bf\ref{tordu}\rm, \textbf{b)}, \V\ is the union of two smooth circles, and of segments
which join these two circles at branch points. And since \V\ can not be decomposed into an union of circle, 
there is no circle
carried by a \vf\ of \V\ which passes over one of these segments (\textbf{c)}), which is a contradiction.

\carr 
\vspace{.3cm}\lt

That is why the existence of a \dix\ prevents the proof to work. That is
also why we have refined the first
cell decomposition of \B\ in subsection \bf\ref{raff}\rm.

\rke For the non compact cells, it is much simpler, since a \tr\ fully
carried by a \vf\ $\R\cI$ can always be
decomposed into a union of smooth lines.
\end{rmk}

\vspace{0.5cm}

\subsection{Resolving sequence of \dec s}\label{suite}

We keep the real number $\eps>0$ defined in subsection \bf\ref{prem}\rm\
for the splitting from \B\ to $\B_1$.
Let $(\eps_n)_{n\in\N}$ be a strictly decreasing sequence of real numbers
such that
for all $n$,
\hbox{$\frac{\eps}{2}<\eps_n<\eps$}.\lt
Let $Y_\eps$ be a $\frac{\eps}{2}$-neighbourood of $Y^{[1]}$ in \B. Then the trace of $\B_1$ on
$\pi^{-1}(Y_\eps)$ is a \sbr\ with boundary denoted $\B'_1$.\lt

The purpose of this subsection is to explain how to build a sequence of \dec s of \B, whose sequence of \dec
s it induces on $\B'_1$ is resolving.\lt

We denote $(y_i)_{i\in J}$ the set of vertices of $Y$, where $J$
is a subset of \N, finite or not. To each vertex $y_i$ correspond several vertices of
$Y_1$, at least 2 and at most 3, wether $y_i$ is a double point or a regular point of the \ls. We denote
these vertices $y_i(j)$ for $j=1$, $2$ or $3$. We then call 
$d_i(j)$ the projection by $p_1$ of the disk
of $\B_1$ centred at $y_i(j)$ and of radius $2\eps$, such that
$d_i=\cup_{j\in \{1,2,3\}}d_i(j)$ is a branched
disk, neighbourhood of $y_i$ in \B\ (see figure \textbf{\ref{didisques}}).

\begin{figure}[!h]
\begin{center}
\input{didisques.pstex_t}
\caption{}
\label{didisques}
\end{center}
\end{figure}

The \ls\ of $\B'_1$ is included in the union of the 
$p_1^{-1}(d_i(j)))$ 's, for all the $i$ 's and $j$ 's.
Notice that the \ls\ of $\B'_1$ has no double point. Moreover, since we supposed that the edges of $Y$ are
strictly longer than $5 \eps$, if $y_{i_1}$ is different from $y_{i_2}$ , then
$d_{i_1}$ and $d_{i_2}$ are disjoint.\\

At last, we define a sequence of vertices of $Y$, $(y_{\psi(n)})_{n\in\N}$, for $\psi$ a
map from \N\ to $J$, such that each vertex appears infinitely many times. This is possible since there is only
countably many vertices in $Y$.\lt

Define now what are the \dec s from $\B_1$ to $\B_2$. We take all the edges of $Y$, a vertex of whose is 
$y_{\psi(1)}$. We orient them from $y_{\psi(1)}$ to their second vertex. Let \ad\ be one of these edges. Its
second vertex is $y_k$, different from $y_{\psi(1)}$.
We call $\V_\ad$ the trace of $\B_1$ on
$\pi^{-1}(\ad)$. Since \ad\ is oriented, it makes sense to talk
 of direct and backward branchings along
$\V_\ad$.

Because of the definition of $\B_1$, the backward 
branchings all lie in $\pi^{-1}(d_{\psi(1)})$, and
the direct branchings all lie in $\pi^{-1}(d_k)$. 
Moreover, each branching lies at a distance \eps\ from
the ends of $\V_\ad$. Actually, at this step of the sequence of \dec s, there is at 
most one direct branching
and one backward branching along $\V_\ad$. If there is no backward branching, no 
\dec\ will be made along 
$\V_\ad$. Else, we will perform a \dec\ along a path inscribed on $\V_\ad$,
going from the backward branching to the direct one if it exists, or to the end of 
$\V_\ad$, in an 
$\eps_1$-neighbourhood of this path. If a direct branching is met, this \dec\ can be an 
over, under or 
neutral \dec. The following subsection \textbf{\ref{adaptes}} will tell which one must be 
chosen. If it is
the neutral \dec, the \dec\ stops at this branching point. Else, 
we can split on along a 
path in $\V_\ad$
which goes to the end of $\V_\ad$. Since $\eps_1 <\eps$, along this path, no other backward 
branching is met,
and hence, there is no backward \dec. The same process is applied to the other edges 
having $y_{\psi(1)}$ as
a vertex. The second \dec\ takes place in an $\eps_2$-neighbourhood of the corresponding 
path, the third \dec\
takes place in an $\eps_3$-neighborhhod of the corresponding path, and so on.

The fact that the $\eps_i$ are decreasing 
allows to avoid backward
\dec s. The order of the edges does not matter.\lt

After these \dec s, we get a \sbr\ $\B_2$. We take all the edges of $Y$, a vertex of whose 
is 
$y_{\psi(2)}$. We orient them from $y_{\psi(2)}$ to their second vertex. Let \ad\ be one of 
these edges. Its
second vertex is $y_k$, different from $y_{\psi(2)}$. We call $\V_\ad$ the trace of 
$\B_2$ on
$\pi^{-1}(\ad)$. The situation is as previously, except for one detail : there can now be 
more than one
direct branching and one backward branching along $\V_\ad$. However, all the backward
branchings lie in 
$\pi^{-1}(d_{\psi(2)})$, and all the direct branchings lie in $\pi^{-1}(d_k)$. All these 
branchings lie at a
distance at least $\eps_i$ from the ends of $\V_\ad$, where $i$ is the number of \dec s 
performed on
$\B_1$.\lt
Look at the backward branchings of $\V_\ad$ : there are $j$ such branchings. Since the 
successive \dec s have
been performed in smaller and smaller neighbourhoods, we can order these branchings from 
the furthest from 
$\pi^{-1}(y_{\psi(2)})$ to the nearest. We note them $b_1\ldots b_j$, $b_i$ being strictly 
further than
$b_{i+1}$. We will make \dec s along paths going from the $b_i$'s, in smaller and smaller 
neighbourhoods,
whose size is set by the $(\eps_n)$ sequence. To avoid any backward \dec, we begin by the 
\dec\ along a path
starting from $b_1$. The second splitting will start from $b_2$, and so on until the last 
\dec, which will
start from $b_j$. When a direct branching is met, one of the over, under and neutral 
\dec s must be chosen
: this is done in subsection \textbf{\ref{adaptes}}. As previously, if the neutral \dec\ 
is chosen, the \dec\
stops here. Else, we can split on until another direct branching is met, or until 
the end of
$\V_\ad$. Again, thanks to the choice of the $\eps_n$, backward \dec s are avoided. 
The same process is
applied to all the edges having $y_{\psi(2)}$ as a vertex. The order of the edges does 
not matter. \lt

\begin{figure}[!h]
\begin{center}
\input{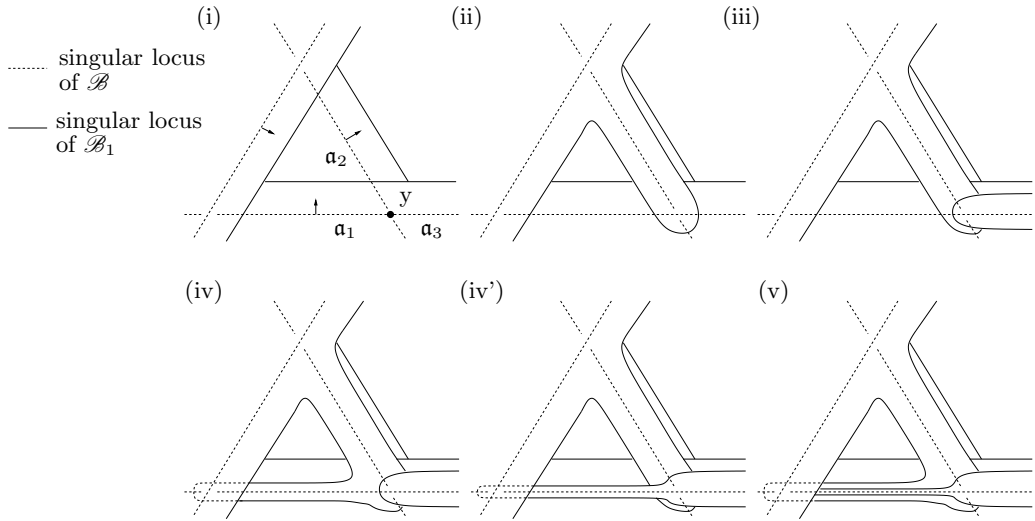}
\caption{Example of a sequence of \dec s}
\label{ancienne}
\end{center}
\end{figure}

We iterate these operations at each step : the \dec s from $\B_n$ 
to $\B_{n+1}$ are 
performed along arcs
whose image by $\pi$ is included in an edge having $y_{\psi(n)}$ 
as a vertex. The backward 
branchings are
always over $d_{\psi(n)}$ : they are more and more numerous, but 
they are always 
strictly ordered, from the
furthest to the nearest. Moreover, the \ls\ of $\B'_{n+1}$ does not 
intersect
$P_{n+1}^{-1}(d_{\psi(n)})$ any more : the singularities over
$d_{\psi(n)}$ have then
been \it resolved\rm .
Since the vertex $y_{\psi(n)}$ will reappear infinitely many times in the
sequence
$(y_{\psi(n)})_{n\in\N}$, the sequence of \dec s is resolving.\lt
Figure \bf\ref{ancienne}\rm\ shows an example of such a sequence of \dec s. On this figure, the branch loci
are seen ``from above", and only the top parts are drawn. 
The three first points show a sequence of \dec s at
the end of which there are several direct branchings along some edges having $y$ as a vertex. The first \dec\
to be performed along $\ad_1$ is the one drawn in \bf(iv)\rm , but not the one drawn in \bf (iv')\rm\
, where a backward \dec\ occurs. The second \dec\ is the one drawn in \bf(v)\rm\ . It then remains to split
along $\ad_2$ and $\ad_3$.

\subsection{Adapted \dec s}\label{adaptes}

We will now see how it is possible to perform the splittings along the edges
previously defined, in such a way that none of the $\B_n$ has a \twc, and then a \dix.\\

If an arc of \dec\ from \B$_n$ to \B$_{n+1}$ is not in a face-to-face
position, then
the \ls\ of \B$_{n+1}$ remains the same as the \ls\ of $\B_n$ : it is
deformed, but
there is no new double point.\lt

When the arc of \dec\ is in a face-to-face position, then we have the
following
fact : an over splitting introduces two new double points in the \ls, a
positive and a
negative one, and an under \dec\ introduces two double points at the same
place
but of opposite signs. Figure
\bf\ref{3dec}\rm\ shows this.
\lt

Being given an arc to split along, we now have to find a \dec\ which
will not
create a \twc. Such a \dec\ will be said \it adapted\rm.\lt

The following proposition is fundamental (we keep the previous notations) :

\begin{prop}\label{fondam} Let $\B_n$ be a \sbr\ obtained from $\B_1$ by a sequence of
\dec s, and which does not have any \twc. We denote
$\Ls_n$ its \ls. Then, for every arc of \dec\ in a face-to-face position in \B$_n'$ (the \sbr\ induced by \B$_n$ on $\pi^{-1}(Y_\eps)$),
at least one
of the three \dec s, over, under or neutral, is adapted.
\end{prop}

\demo According to corollary \textbf{\ref{simple}}, it is possible to consider only simple \twc s. 
This will be done throughout this proof.\lt
Let us denote \ad\ the \dec\ path used to go from $\B_n$ to $\B_{n+1}$. We suppose that \ad\ lies in a face-to-face situation. Consider all the possible \dec s along \ad. When we perform one of these \dec s, the \ls\ is only modified in some neighbourhood \vad\ of \ad.\lt

\begin{figure}[!h]\begin{center}\input{nota.pstex_t}\caption{}\label{nota}
\end{center}
\end{figure}

We begin by performing the neutral \dec. We note \Bn\ the obtained \sbr. We use the notations of figure \textbf{\ref{nota}}, \textbf{b)}, where $\ad_g=[p_g,q_g]$, $\ad_d=[q_d,p_d]$, $l_g$ and $l_d$ are smooth oriented segments of the \ls. We will say that an oriented curve passes \textit{positively} through one of these segments if it passes through this segment with the same orientation than this segment. We will say that it passes \textit{negatively} if it passes through a segment in the other way.\lt

Suppose that the neutral \dec\ is not adapted. Suppose, for instance, that \ga\ is a positive \twc\ passing positively through $\ad_g$ in \Bn. The other cases can be dealt with in the same way. We distinguish the case \textbf{(A)}, where \ga\ does not pass through $\ad_d$, from the case \textbf{(B)}, where it passes through $\ad_d$.\lt

\begin{itemize}
\item[]
\begin{leme}\label{leme1}
If $\gamma$ is in case \textbf{(B)}, then it passes negatively through $\ad_d$.
\end{leme}
\vspace{.2cm}
\noindent \bf Proof \rm 
Suppose that it passes positively through $\ad_d$. The immersion of \ga\ in \Bn\ gives an immersion $\ga'$ of a closed curve in $\B_n$, which coincides with \ga\ outside \vad. It remains to define $\ga'$ inside \vad. Let us define two oriented edges in $\B_n$, denoted $a_g$ and $a_d$, which go from $p_g$ to $q_g$ and from $q_d$ to $p_d$.
Those edges correspond to $\ad_g$ and $\ad_d$, but they are not in the \ls\ of $\B_n$. We define $\ga'$ inside \vad\ in the same way as \ga, replacing $\ad_g$ and $\ad_d$ by $a_g$ and $a_d$. Each corner of $\ga'$ different from $p_g$, $p_d$, $q_g$ and $q_d$ is ascending. If we consider $\ga'$ as a loop based in $q_g$, we can write $\gamma'=\gamma'_1\ast a_d\ast\gamma'_2\ast a_g$, where $\gamma'_1$ is the part of$\gamma'$ going from $q_g$ to $q_d$, and $\gamma'_2$ is the part of $\ga'$ going from $p_d$ to $p_g$. We then set :
$\beta_1=\gamma'_1\ast [q_d,q_g] $ and $\beta_2=\gamma'_2\ast [p_g,p_d]$. These two curves lie in the \ls\ of 
$\B_n$. Moreover, a corner of $\ga'$ is necessarily a corner of $\beta_1$ or of $\beta_2$, and \it 
vice-versa\rm . Therefore, one of the curves $\beta_1$ or $\beta_2$ has at least one corner, and all of its corners are ascending. It is a positive \twc\ of $\B_n$, existing before \dec, which is a contradiction.\carr
\end{itemize}

 \vspace{.4cm}

In case \textbf{(B)}, we can write $\gamma=\gamma_1\ast\ad_d^{-1}\ast\gamma_2\ast\ad_g$, where 
$\gamma_1$ is the part of $\gamma$ going from $q_g$ to $p_d$, and $\gamma_2$ the one from $q_d$ to $p_g$. We have a second lemma :
 \vspace{.4cm}

\begin{itemize}

\item[]

\begin{leme}\label{leme2}
There is at least one corner on $\ga_1$ and at least one corner on $\gamma_2$.
\end{leme}

 \vspace{.2cm}

\noindent \bf Proof \rm 

Suppose for instance that there is no corner on $\ga_2$. As in the proof of lemma \textbf{\ref{leme1}}, there is a closed curve $\ga'$ in $\B_n$ corresponding to \ga, and which can be written $\gamma'=\gamma'_1\ast a_d^{-1}\ast\gamma'_2\ast a_g$. We set in $\B_n$ : $\beta=\gamma'_1\ast [p_d,p_g]\ast\gamma_2^{'-1}\ast [q_d,q_g]$. This is a closed curve, immersed in $\B_n$, and included in the \ls\ of $\B_n$. Each corner of \ga\ is a corner of $\ga'_1$, and thus of $\beta$, and each corner of $\beta$ is a corner of \ga. Hence, $\beta$ is a positive \twc\ in $\B_n$, which is a contradiction.\carr

\end{itemize}

  \vspace{.4cm}

\begin{figure}[!h]\begin{center}\input{case.pstex_t}\caption{Positive \twc s passing through \vad}
\label{case}
\end{center}
\end{figure}

In any case, we will prove that the over \dec\ is adapted. We perform this \dec. We call \Bs\ the obtained \sbr.
We use the notations of figure \textbf{\ref{nota}}, \textbf{c)}, where $q_g$ and $q_d$ are double points of the \ls\ of \Bs. Suppose that \Bs\ contains a positive simple \twc\ $\delta$ passing inside \vad, i.e. the over \dec\ makes $\delta$ appear. We will actually prove that this cannot happen. Figure 
\textbf{\ref{case}} shows all the possible local configurations of positive simple \twc s passing in \vad. If we revert the orientations of these curves, we get all the possible local configurations of negative simple \twc s passing in \vad.\lt

Diagrams \textbf{1} and \textbf{3} of this figure are equivalent, in the sense that there is a positive simple \twc\ as in 
\textbf{1} \ssi\ there is a positive simple \twc\ as in \textbf{3}. In the same way, diagrams \textbf{5} and \textbf{7} are equivalent, as diagrams \textbf{2} and \textbf{4} and diagrams \textbf{6} and \textbf{8}. There are thus only 4 cases to study for $\delta$.\lt
Notice that the immersion of $\ga$ into \Bn\ is also an immersion of \ga\ into \Bs. To avoid confusions, the image of this last immersion is called \gs. However, $q_g$ and $q_d$ are smooth points of \ga, whereas they are descending corners of \gs. The curve \gs\ is not twisted, but all its corners different from $q_g$ and $q_d$ are ascending.\lt
We now study all the possible cases, starting with those where \ga\ is in case 
\hbox{\textbf{(A)} :}
\vspace{.5cm}
\begin{itemize} 
\item[(A.1)] 
\end{itemize}
\begin{figure}[!h]\begin{center}\input{A1.pstex_t}\caption{(A.1)}\label{A1}
\end{center}
\end{figure}

Figure \textbf{\ref{A1}} shows what happens. We consider \gs\ and $\delta$ as two loops based in $q_g$, and we set 
$\beta=\delta\ast\gs$. This loop is freely homotopic in \Bs\ to an immersed loop which does not pass in \vad\ anymore, and whose corners are all ascending. This last loop contains, according to lemma \textbf{\ref{simple}}, a positive simple \twc\ which does not pass in \vad\ neither. This curve is a positive simple \twc\ in $\B_n$, which is a contradiction.

\vspace{.5cm}

\begin{itemize} 
\item[(A.4)] 
\end{itemize}
\begin{figure}[!h]\begin{center}\input{A4.pstex_t}\caption{(A.4)}\label{A4}
\end{center}
\end{figure}

Figure \textbf{\ref{A4}} shows what happens. We consider \gs\ and $\delta$ as two loops based in $q_g$, and we set 
$\beta=\delta\ast\gs$. This loop is immersed, and $q_g$ is no longer a corner of $\beta$. All the corners of $\beta$ are ascending. This $\beta$ contains a positive simple \twc. If this curve passes in \vad, it passes in a row either through  $\ad_g$, $l_h$ and $\ad_d$, positively, or through $l_d$, $l_b$ and $l_g$ positively. This curve is a positive simple \twc\ in $\B_n$, which is a contradiction.\vspace{.5cm}\lt

\begin{itemize} 
\item[(A.6)] 
\end{itemize}
\begin{figure}[!h]\begin{center}\input{A6.pstex_t}\caption{(A.6)}\label{A6}
\end{center}
\end{figure}

Figure \textbf{\ref{A6}} shows what happens. We consider \gs\ and $\delta$ as two loops based in $q_g$, and we set 
$\beta=\delta\ast\gs$. This $\beta$ is homotopic in \Bs\ to an immersed loop for which $q_g$ is not a corner, and which 
contains a positive simple \twc. If this curve passes in \vad, it passes in a row either through $l_d$, $l_b$ and $l_g$ positively, and nowhere else in \vad. This curve is a positive simple \twc\ in $\B_n$, which is a contradiction.

\vspace{.5cm}
\begin{itemize} 
\item[(A.7)] 
\end{itemize}
\begin{figure}[!h]\begin{center}\input{A7.pstex_t}\caption{(A.7)}\label{A7}
\end{center}
\end{figure}

Figure \textbf{\ref{A7}} shows what happens. We consider \gs\ and $\delta$ as two loops based in $q_g$, and we set 
$\beta=\delta\ast\gs$. This $\beta$ is homotopic in \Bs\ to an immersed loop for which $q_g$ is not a corner, and which 
contains a positive simple \twc. If this curve passes in \vad, it passes in a row either through $\ad_g$, $l_h$ and $\ad_d$ positively, and nowhere else in \vad. This curve is a positive simple \twc\ in $\B_n$, which is a contradiction.

\vspace{.5cm}

We then deal with the cases where \ga\ is in case \textbf{(B)} :
\vspace{.5cm}
\begin{itemize} 
\item[(B.1)] 
\end{itemize}

In the same way as we can build \gs\ from \ga, we can build a closed curve \dn\ immersed in \Bn, from $\delta$. This curve has the same corners as $\delta$, minus $q_g$. If $q_g$ is the only corner of $\delta$, then \dn\ has no corner, else it is a positive \twc.\lt
We consider \dn\ and $\ga$ as two loops in \Bn\ based in $q_g$, and we set $\beta=\dn\ast\ga$.
This loop is freely homotopic in \Bs\ to an immersed loop which does not pass through $\ad_g$, and whose corners are all ascending. This last loop contains, according to lemma \textbf{\ref{simple}}, a positive simple \twc\ which does not pass through $\ad_g$ neither. We still call this curve $\beta$. If this curve does not pass in \vad, it is a positive simple \twc\ in $\B_n$, which is a contradiction.\lt
If this curve passes in \vad, it passes in a row either through $l_d$ and $\ad_d$, positively, and nowhere else in \vad. The immersion of $\beta$ in \Bn\ implies the existence of the immersion of a closed curve $\beta_{sup}$ in \Bs, and which passes only through $\ad_d$ and $l_d$ in \vad. 
We then modify $\delta$ by adding to it a loop $l_b\ast l_h$, so that we get a curve $\delta'$ modeled on diagram \textbf{4} of figure \textbf{\ref{case}}. \lt
We get a contradiction in the same manner as in point 
\textbf{A.4}, by using $\beta_{sup}$ and $\delta'$.

\vspace{.5cm}
\begin{itemize} 
\item[(B.2)] 
\end{itemize}
By symmetry, this point can be dealt in the same way as the previous one, \textbf{B.1}.

\vspace{.5cm}
\begin{itemize} 
\item[(B.7)] 
\end{itemize}
\begin{figure}[!h]\begin{center}\input{B7.pstex_t}\caption{(B.7)}\label{B7}
\end{center}
\end{figure}

This case is shown on figure \textbf{\ref{B7}}. As previously, we write $\gs=\ga_1\cup\ga_2$, where $\ga_1$ has $q_g$ as first end and $q_d$ as last end. In \Bs, we have the following loop, based in $q_g$ :$\beta=\delta\ast l_b\ast\ga_2$. Neither $q_g$ nor $q_d$ are corners of this loop. The corners of this loop are the ones of $\ga_2$ and the ones of $\delta$, except $q_d$. They are thus all ascending, and according to lemma \textbf{\ref{leme2}}, there is at least one. Thus,  $\beta$ is a positive simple \twc, which was there before \dec, which is a contradiction.

\vspace{.5cm}
\begin{itemize} 
\item[(B.8)] 
\end{itemize}
By symmetry, this point can be dealt in the same way as the previous one, \textbf{B.7}.

\vspace{.5cm}

After these 8 points, the over \dec\ is adpated.\lt
The other cases where the neutral \dec\ is not adapted are dealt with in the same way, and are the following :

\begin{itemize}
\item[-] $\gamma$ is positive and passes negatively through $\ad_g$ : the under \dec\ is adapted ;
\item[-] $\gamma$ positive and passes positively through $\ad_d$ : the under \dec\ is adapted ;
\item[-] $\gamma$ is positive and passes negatively through $\ad_d$ : the over \dec\ is adapted.
\end{itemize}
\carr

\subsection{Conclusion}

After the two previous subsections, we have built a sequence of \dec s of
\B, none
of whose having a \dix. This sequence induces a resolving sequence of \dec
s of
$\B'_1$, whose inverse limit is a lamination $\lambda$ fully caried by
$\B'_1$. We
aim at proving that $\lambda$ has null holonomy.\lt

Let $\Sigma$ be a 2-cell of $Y$, and $\partial\Sigma\cI$ be the subbundle
of \NB\
over $\partial\Sigma$. Then $\lambda\cap(\partial\Sigma\times I)$ is an
oriented
dimension 1 lamination denoted $l_\Sigma$, fully carried by $\partial\Sigma\times I$, and
obtained as
the inverse limit (in the sense of definition \bf\ref{liminv}\rm ) 
of the oriented \tr s
$v_n=\B_n\cap(\partial\Sigma\times I)$.\\
If $\Sigma$ is not compact, there are no holonomy problems since there is no
first-return map on a fibre. \lt
Then, suppose that $\Sigma$ is compact.

\defi Let $\lambda$ be an oriented lamination carried by a trivial bundle $\SU\cI$. An \it increasing leaf
\rm (resp. a \it decreasing leaf\rm ) of $\lambda$ is a leaf which goes, in the direct way, from a point 
$p_1=(\theta,t_1)$ to a point $p_2=(\theta,t_2)$, with $t_1<t_2$ (resp. $t_1>t_2$).
\end{df}

\begin{leme} The lamination $l_\Sigma$ is a lamination by circles.
\end{leme}

\demo We denote $N(v_n)=N(\B_n)\cap(\partial\Sigma\times I)$, which is a \vf\ of $v_n$. 
We call $L_\Sigma=\cap_{n\in\N}N(v_n)$, and we then have
$l_\Sigma=\partial L_\Sigma$, according to definition \bf\ref{liminv}\rm . \lt
Let $L$ be an increasing leaf of $l_\Sigma$. This leaf is a spiral with two limit circles :
$C^+$, limit when
$L$ is followed in the direct way, and $C^-$, limit when $L$ is followed in the indirect 
way. We call $A$ the
annulus between $C^+$ and $C^-$. Look at 
$L_\Sigma\cap A$. By construction, this intersection is notequal to $A$. This means that 
$A\backslash L_\Sigma$ contains some subset of the form $\gamma\cI$, where $\gamma$ is a 
compact oriented
path fully carried by $\partial\Sigma\times I$, and which is increasing 
(see. figure \textbf{\ref{tour}}). 

\begin{figure}[!h]\begin{center}\input{tour.pstex_t}\caption{}\label{tour}
\end{center}\end{figure}

Hence, there exists an integer $N$ such that for all integer $n$ greater than $N$, we have 
$(\gamma\cI )\cap N(\B_n)=\emptyset$. 
If not, it
would imply the existence of a sequence of points $(q_n)$ such that $q_n\in 
(\gamma\cI )\cap N(\B_n)$. Since
$(\gamma\cI )\cap N(\B_n)$ is compact, there would be a subsequence of $(q_n)$ converging 
towards a point $q$
contained in $\cap_{n\in\N}((\gamma\cI )\cap N(\B_n))$. But this last set is equal to 
$\gamma\cap L_\Sigma$, which
is empty. \lt
It is then impossible to find a path in $v_N$ going in the direct way from $p_2$ to $p_1$, where $p_1$ and
$p_2$ are two points of $L$ placed as in figure 
\textbf{\ref{tour}}.\lt
However, because $\Sigma$ is not a \dix\ and according to lemma \textbf{\ref{voi1}}, the existence of a path 
of $v_N$ going in the direct way from $p_1$ to $p_2$ implies the existence of a path 
of $v_N$ going in the indirect way from $p_1$ to $p_2$. This is a contradiction, and $L$ must be a circle.
In the same way, $l_\Sigma$ does not have any decreasing leaf.
\carr
\vspace{.2cm}

Hence, $\lambda$ has null holonomy.\lt

To get a lamination fully carried by \B, it only remains to ``fill the
holes" of leaves
of $\lambda$, these holes being in fact diffeomorphic to the 2-cells of
$Y$, which
are disks and half planes. This is possible because $\lambda$ is null
holonomic. This ends the proof of theorem \bf\ref{but}\rm .

\section{Some remarks about the problem of determining whether a \sbr\ fully carries a lamination or not}

Proposition \textbf{\ref{nonport}} states that a necessary condition for a \sbr\ to fully carry a lamination is that it does not have any \dix\ which is a sector. The obstruction to the existence of a lamination fully carried seems to be essentially linked to the phenomenon of the twisted disks of contact. However, we will show on two examples that the non-existence of twisted disks of contact is not a necessary condition nor a sufficient one for a \sbr\ to fully carry a lamination.

\subsection{No \dix\ is not necessary }

Let \B\ fully carry a lamination $\lambda$ and have a \dix\ \di. According to proposition
\textbf{\ref{nonport}}, \di\ cannot be a sector. Then using the arguments of the proof of 
proposition \textbf{\ref{nonport}}, over the boundary of \di , the trace of some leaf $l$ of $\lambda$ must 
be a spiral. Hence, $l$ cannot pass all over \di. That means that $l$ passes over \di\ along some annulus 
around $\partial\di$, and then quits \di. This separation of $l$ and \di\ can be done only along a piecewise 
smooth circle \SB\ of the branch locus such that \SB\ and $\partial\di$ bound a sink annulus in \di. 
Moreover, if we look at all the connected surfaces immersed in \B\ and bounded by \SB\ so that the branch 
directions along the boundary component \SB\ point outwards, then one of these surfaces is not a disk. If 
they were all disks, $l$ could not exist, again with the arguments of the proof of proposition 
\textbf{\ref{nonport}}. This surface which is not a disk allows the existence of the holonomy of the spiral 
traced by $l$ over $\partial\di$, because its $\pi_1$ is not zero. The \sbr\ of figure \textbf{\ref{notnecess}} is an example of such a 
\sbr.

\begin{figure}[!h]\begin{center}\input{notnecess.pstex_t}\caption{A \sbr\ with a \dix\ fully carrying a
lamination}\label{notnecess}
\end{center}
\end{figure}

\subsection{No \dix\ is not sufficient}

Indeed, a \dix\ may be ``hidden'' in the \sbr , as shown in the example of figure \textbf{\ref{dcvv}}.\lt

\begin{figure}[!h]\begin{center}\input{dcvv.pstex_t}\caption{}\label{dcvv}
\end{center}
\end{figure}

\begin{figure}[!h]\begin{center}\input{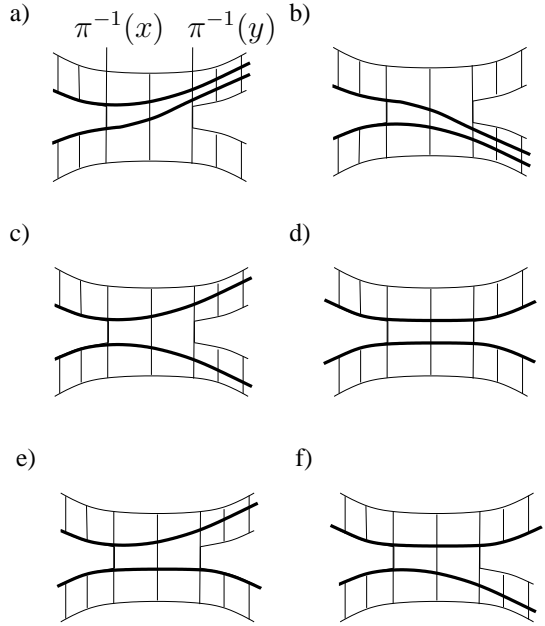}\caption{Trace of two boundary leaves of $\lambda$}\label{passe}
\end{center}
\end{figure}

This figure shows the \ls\ of a \sbr. This \sbr\ cannot fully carry a lamination, although it has no \dix. Indeed, were it the case, a lamination $\lambda$ fully carried would have to boundary leaves (i.e. there is no other leaf between them and some component of the horizontal boundary) passing through $\pi^{-1}([x,y])$, as in one of the situations of figure  \textbf{\ref{passe}}.

In case \textbf{a)} of this figure, the \sbr\ obtained after an over \dec\ along $[x,y]$ still fully carries $\lambda$. In case
\textbf{b)} only the under \dec\ has this property, in case \textbf{c)} the three \dec s fit, in case \textbf{d)} only the neutral \dec\ fits, in case \textbf{e)} the neutral and the over \dec s fit, and in case \textbf{f)} the neutral and the under \dec s fit. However, if any of these \dec s is performed along $[x,y]$, a \dix\ appears and the new \sbr\ cannot fully carry a lamination.\lt

Actually, it is possible to give a definition of this kind of \suct, which always give birth to a \dix, whatever is the \dec\ performed along some arcs. The non existence of this kind of \suct\ could thus be a sufficient condition for a lamination fully carried to exist. More generally, it is natural to ask whether a \sbr\ without any \suct\ fully carries a lamination or not. We could then try to adapt the proof of theorem \textbf{\ref{but}}, that is to perform an infinite sequence of \dec s called adapted, which does not create any \dix, or even any \suct. It is possible to show that a \dec\ along an arc \ad\ cannot create a \suct\ on one side of \ad\ and another one of the same sign of the other side of \ad. But I cannot find an argument to show that there exists a \dec\ which does not create a \suct\ on each side, but also which does not create two \sucts\ of the same sign and on the same side of \ad. The issue is that the boundary of a \suct\ can be too ``wild", and can pass several times through \ad. A way to avoid this problem is to consider \twc s, bounding or not a surface, and to use lemma \textbf{\ref{simplet}} and corollary \textbf{\ref{simple}}.

\end{document}